\documentclass[11pt]{article}
\usepackage{a4wide}
  
\title{There are only two nonobtuse binary triangulations of the unit $n$-cube}
\date{}
\author{Jan Brandts, Sander Dijkhuis, Vincent de Haan, and Michal K\v{r}\'{\i}\v{z}ek}

\usepackage{epsfig}
\usepackage{amsmath}
\usepackage{amssymb} 
\usepackage{makeidx}
\usepackage{graphics}
\usepackage{color} 
 
\newtheorem{Th}{Theorem}[section]
\newtheorem{Le}[Th]{Lemma}
\newtheorem{Co}[Th]{Corollary}

\newtheorem{Def}[Th]{Definition}
\newtheorem{rem}[Th]{Remark}

\newcommand{\be}{\begin{equation}}
\newcommand{\ee}{\end{equation}}
\newcommand{\RR}{\mathbb{R}}
\newcommand{\sth}{\hspace{2mm} | \hspace{2mm}}
\newcommand{\half}{\frac{1}{2}}
\newcommand{\hdrie}{\hspace{3mm}}
\newcommand{\und}{\hdrie\mbox{\rm and }\hdrie}

\begin{document}
 
\maketitle

\def\BB {\mathcal{B}}
\def\SS {\mathcal{S}}
\def\NN {\mathcal{N}}
\def\TT {\mathcal{T}}
\def\OO {\mathcal{O}}
\def\conv {\mbox{\rm conv}}
\newcommand{\supp}{{\rm supp}}
\newcommand{\ol}{\overline}

\vspace*{-10mm}
\begin{abstract}
Triangulations of the cube into a minimal number of simplices without additional vertices have been studied by several authors over the past decades. For $3\leq n\leq 7$ this so-called simplexity of the unit cube $I^n$ is now known to be $5,16,67,308,1493$, respectively. In this paper, we study triangulations of $I^n$ with simplices that only have nonobtuse dihedral angles. A trivial example is the standard triangulation into $n!$ simplices. In this paper we show that, surprisingly, for each $n\geq 3$ there is essentially only one other nonobtuse triangulation of $I^n$, and give its explicit construction. The number of nonobtuse simplices in this triangulation is equal to the smallest integer larger than $n!({\rm e}-2)$.
\end{abstract}
{\small{\bf Keywords:} triangulation, simplexity, unit cube, binary simplex, acute simplex, nonobtuse simplex, orthogonal simplex, path simplex, binary matrix, $0-1$ matrix, $0/1$ polytope.\\[1mm]
{\bf 2000 Mathematics Subject Classification:} 52C22, 52B05, 05A15.}
  
\section{Introduction}
Since the paper \cite{Mar} by Mara there has been an interest in finding minimal triangulations of the unit $n$-cube $I^n\subset\RR^n$, where $I=[0,1]$. With triangulations we mean face-to-face partitions of $I^n$ into a minimal number of nondegenerate $n$-simplices, whose set of $n+1$ vertices is a subset of the set of $2^n$ vertices of $I^n$. We will call such simplices binary and denote the set of all binary $n$-simplices by $\BB^n$. The cardinality $\beta_n$ of $\BB^n$ is obviously bounded by the number $\sigma_n$ of all subsets of $n+1$ elements of $\{1,\dots,2^n\}$, and this bound is asymptotically attained \cite{Kom} in the sense that
\be \lim_{n\rightarrow\infty} \frac{\beta_n}{\sigma_n} = 1. \ee
Thus, in theory, for any fixed value of $n$, all binary triangulations of $I^n$ can be found by brute force in finite time. However, both $\sigma_n$ and $\beta_n$ grow very rapidly, as shown in the table below. The values for $\beta_n$ in Table 1, or rather the differences $\sigma_n-\beta_n$, were taken from \cite{Rak}.
\[ \begin{array}{|c||r|r|r|}
\hline
n & \sigma_n & \beta_n & \nu_n\\
\hline
\hline
1 & 1 & 1 & 1\\
\hline
2 & 4 & 4 & 4\\
\hline
3 & 70 & 58 & 34\\
\hline
4 & 4368 & 3008 & 480\\
\hline
5 & 906192 &  556192 & 9984\\
\hline
6 & 621216192 & 366179200 &  284672 \\
\hline
7 & 1429702652400 & 858240222176 &  10474336\\
\hline
\end{array}
\]
\begin{center}
{\bf Table 1. } All, all nondegenerate, and all nondegenerate nonobtuse binary simplices.
\end{center}
Construction of binary cube triangulations turns out to be a highly non-trivial mathematical problem with for example relations to the maximum Hadamard determinant problem \cite{Hada}. This is why the cardinalities $\tau(n)$ of minimal triangulations of $I^n$ have so far only been found \cite{Cot,CoBr,Hug,HuAn} for dimensions $n\leq 7$. For higher dimensions, only bounds and estimates are available, see for instance \cite{Gla,OrSa,Smi}.
\[
\begin{array}{|r||r|r|r|r|r|r|r|}
\hline
n & 1 & 2 & 3 & 4 & 5 & 6 & 7 \\
\hline
\hline
\tau(n) & \phantom{142}1 & \phantom{142}2 & \phantom{142}5 & \phantom{14}16 & \phantom{14}67 & \phantom{1}308 & 1493\\
\hline
\end{array}
\]
\begin{center}
{\bf Table 2. } Cardinality $\tau(n)$ of the minimal binary triangulation of $I^n$.
\end{center}
In this paper we will triangulate the $n$-cube using only nonobtuse binary simplices. These are binary simplices whose dihedral angles are not obtuse, and thus either acute or right. This is illustrated for $n=3$ in Figure 1. On the left we see a binary tetrahedron with only acute angles between its facets. It is, in fact, regular. The one on the right has an obtuse angle between its bottom and front facet. The remaining two have both acute and right angles. Therefore, the left three tetrahedra are nonobtuse, whereas the right one is not.
\begin{center}
\includegraphics[height=3cm]{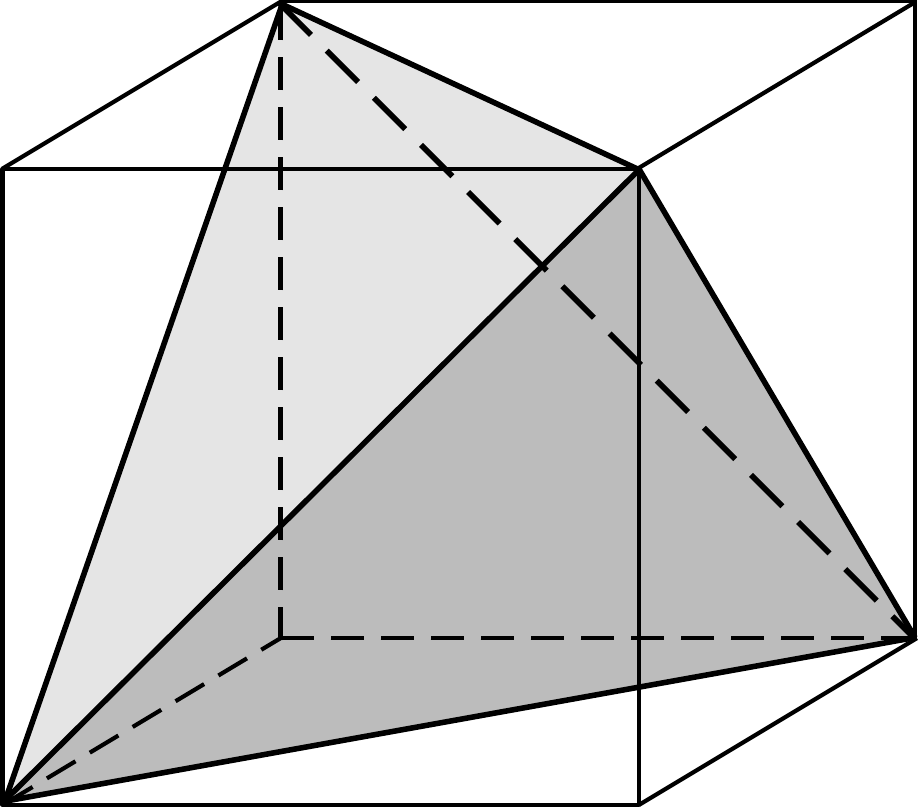}\hspace*{4mm}\includegraphics[height=3cm]{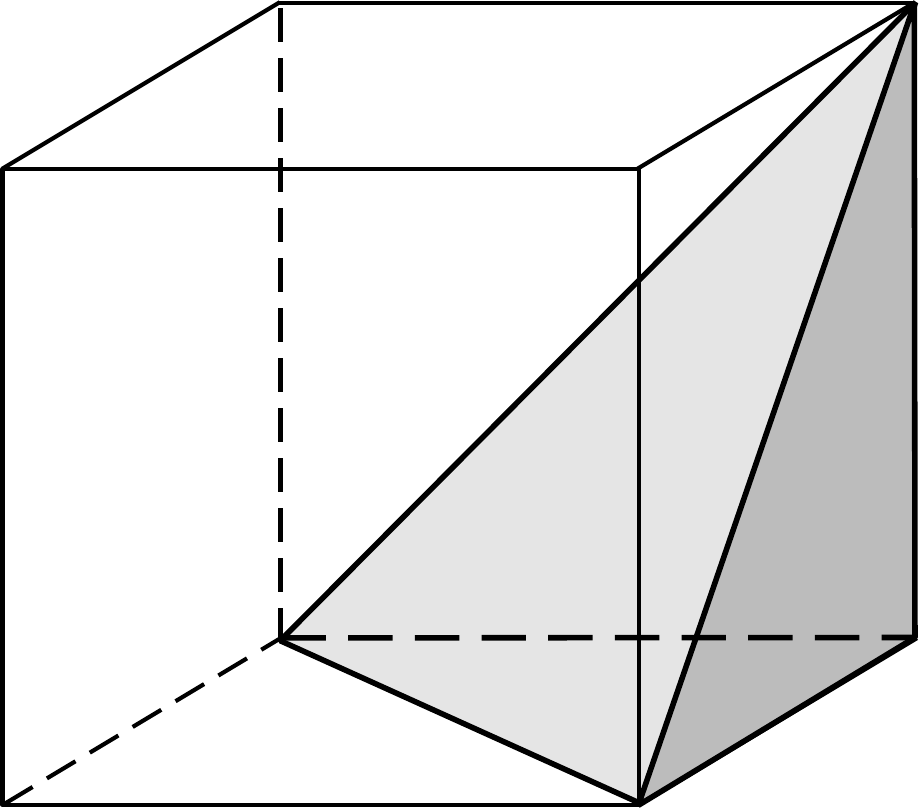}\hspace*{4mm}\includegraphics[height=3cm]{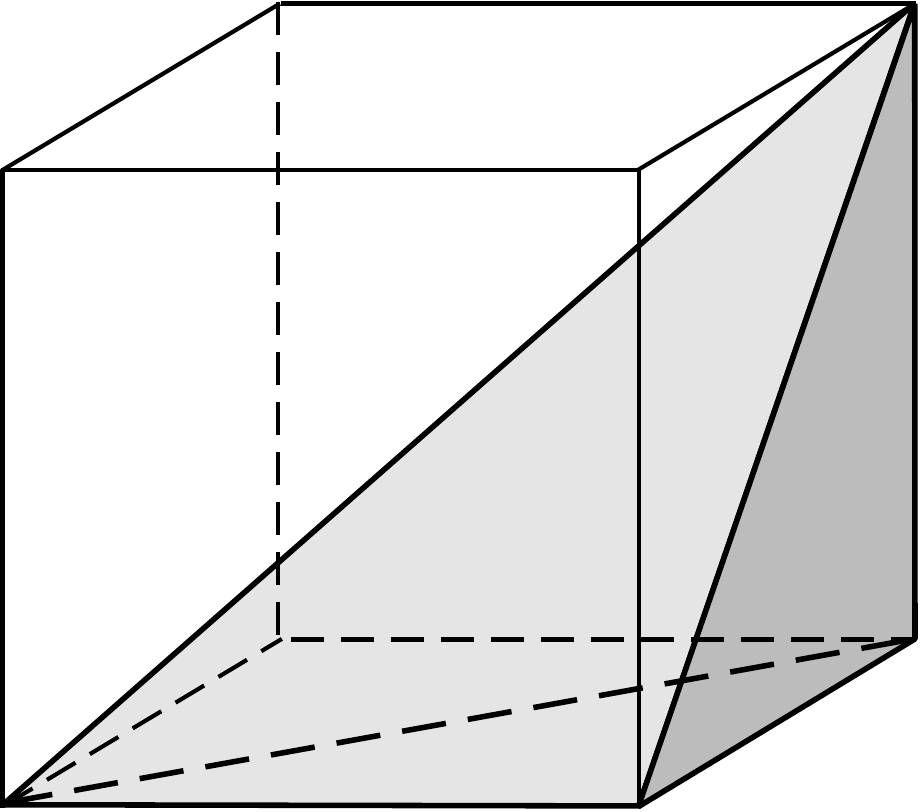}\hspace*{4mm}\includegraphics[height=3cm]{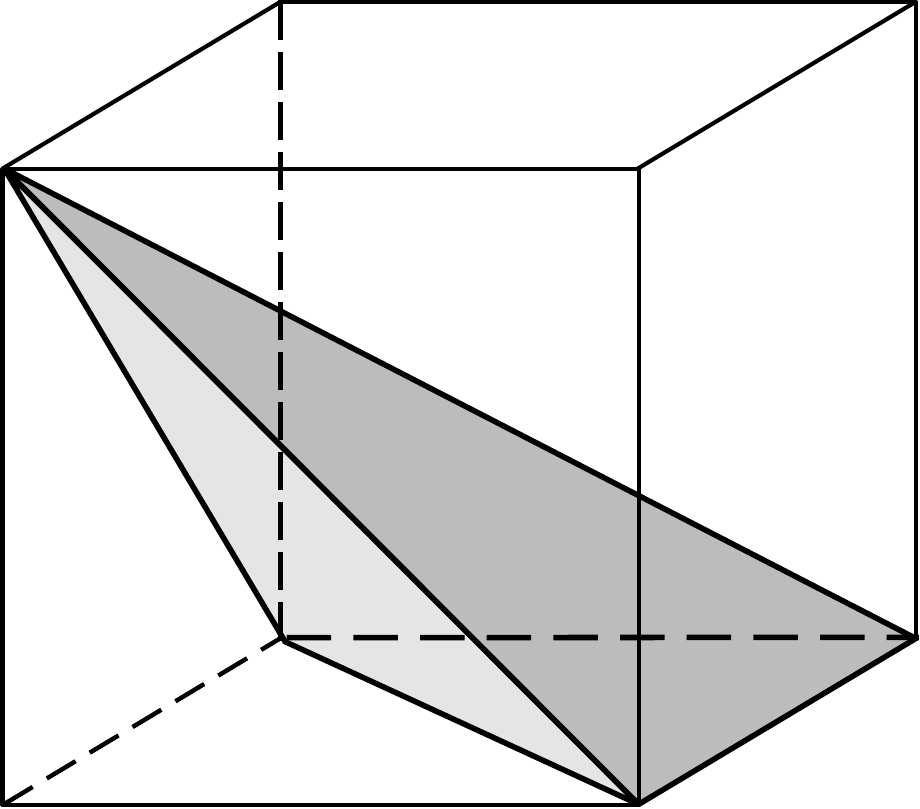}\\
{\bf Figure 1.} Three nonobtuse and one obtuse binary tetrahedron.
\end{center}
For $n\leq 7$, the numbers $\nu_n$ of nonobtuse binary $n$-simplices are listed in the right column of Table 1. They were computed by the authors using exhaustive enumeration, whose details are outside the scope of this paper. What matters is that $\nu_n$ is drastically smaller than $\beta_n$. Even more important is that the restriction of nonobtuseness considerably limits the number of feasible neighbors of a simplex in a triangulation. As a consequence, we were able to completely solve the nonobtuse binary triangulation problem in any dimension.\\[2mm]
For brevity, we will occasionally omit the adjective ``binary'' in the context of binary simplices, assuming that no confusion will arise.\\[2mm]
{\bf Summary of the main results.} {\sl Let $n\geq 3$. Apart from the standard triangulation $\TT_\SS^n$ of $I^n$ into $n!$ path-simplices, there exists another family of nonobtuse binary triangulations $\TT_\NN^n$ of $I^n$. It consists of $\NN(n)$ simplices, where
\be \NN(n) = n\NN(n-1)- n+2 \hdrie \mbox{\rm with }\hdrie \NN(3)=5. \ee
Conversely, each nonobtuse binary triangulation of $I^n$ is essentially equal to $\TT_\SS^n$ or $\TT_\NN^n$. Consequently, $\TT_\NN^n$ is the minimal nonobtuse binary triangulation of $I^n$.
The first few values of $\NN(n)$ are displayed in Table 3.
\[
\begin{array}{|r||r|r|r|r|r|r|r|}
\hline
n & 3 & 4 & 5 & 6 & 7 & 8 & 9\\
\hline
\hline
\NN(n) & \phantom{142}5 & \phantom{14}18 & \phantom{14}87 & \phantom{1}518 & 3621 & 28962 & 260651\\
\hline
\end{array}
\]
\begin{center}
{\bf Table 3. } Cardinality $\NN(n)$ of the minimal nonobtuse binary triangulation of $I^n$.
\end{center}
Alternatively, $\NN(n)$ can be expressed explicitly as
\be \NN(n) = 1+ n!\sum_{k=2}^{n} \frac{1}{k!}, \und \lim_{n\rightarrow\infty} \frac{\NN(n)}{n!} = {\rm e}-2 = 0.718281\dots\ee
In fact, it is easy to show that $\NN(n)$ is equal to the smallest integer larger than $n!({\rm e}-2)$}.\\[2mm]
The triangulation $\TT_\NN^3$ is of course well known, and consists of the regular simplex in the left picture in Figure 1 together with its four neighboring corners of the cube. But as far as we know, its higher dimensional relatives $\TT_\NN^n$ for $n\geq 4$ have not yet been explicitly recognized in the literature as being nonobtuse.

\subsection{Related topics}
A related topic of interest is the triangulation of the cube in simplices having only acute dihedral angles, which introduces even more severe restrictions. It can be done in dimensions two \cite{Gar} and three \cite{KoPaPr,VaHiZhGu} if and only if one allows additional vertices to be introduced. In dimensions four and up, there exists no acute triangulation of $I^n$. Considering the unit cube as a particular instance of a $0/1$-polytope, another related topic would be the nonobtuse or acute triangulation of $0/1$-polytopes in general. One could ask the question which $0/1$ polytopes can be triangulated with nonobtuse simplices, and in how many different ways this can be done.  We refer to \cite{Zie} for a good introduction to $0/1$-polyoptes and to the recent book \cite{San} for triangulation issues in general.

\subsection{Outline of this paper}
After recalling some basic properties of simplices in Section \ref{sect2.1} we specialize the results to nonobtuse binary simplices in Section \ref{sect2.2}. Particular attention will be paid to orthogonal binary simplices. These have a spanning tree of mutually orthogonal edges that form a subset of the edges of $I^n$. Each leaf of the tree turns out to be a vertex opposite an exterior facet, which is a facet of the simplex that lies in a facet of $I^n$. Well-known orthogonal binary simplices are cube corners and path simplices, of which examples are depicted in the second and the third picture of Figure 1, respectively. Also in this paper they will play a crucial role in the analysis. We will encounter them in Section \ref{sect2.3}, together with their fake counterparts: simplices that look like a cube corner or a path simplex from the exterior of $I^n$, but which are not. In Section \ref{sect3} we construct a nonobtuse triangulation of $I^n$ that for $n\geq 4$ we did not encounter in the literature. The hardest part of this paper is Section \ref{sect4}, where we show that each nonobtuse triangulation of $I^n$ must contain either a path simplex, or a cube corner. Moreover, we prove that once a path simplex or a cube corner in $I^n$ is fixed, there is a unique way to finish the triangulation. In case of the cube corner, this leads to our new family of triangulations $\TT_\NN^n$, and in case of the path simplex, it leads to the standard triangulation $\TT_\SS^n$ into $n!$ path simplices. 

\section{Nonobtuse binary simplices}
Here we study nonobtuse binary simplices. In particular, we introduce orthogonal binary simplices, such as path simplices and cube corners. They also play an important role in the general binary simplex triangulation problem. See for instance \cite{CoBr,Cot,Hug,HuAn}. Also in numerical analysis, combinatorics, and computational geometry, these special simplices appear. See for instance \cite{BrKoKrSo} for an overview.

\subsection{Simplices in general}\label{sect2.1}
We recall some basic results on general simplices. In Sections \ref{sect2.2} and \ref{sect2.3} we will specialize these results to binary simplices. We refer to Fiedler \cite{Fie1, Fie3, Fie2} for the original proofs based on graph theory. Linear algebraic proofs can be found in \cite{BrKoKr}.

\begin{Le}\label{pro1}{\sl Each $n$-simplex $S$ has at least $n$ acute dihedral angles. Moreover, if $S$ is nonobtuse, then:\\[2mm]
$\bullet$ each facet of $S$ is a nonobtuse $(n-1)$-simplex;\\[2mm]
$\bullet$ each vertex $p$ of $S$ projects orthogonally into its opposite facet $F$.\\[2mm]
The converse implication holds only for the second statement.\hfill $\Box$}
\end{Le}
In the above lemma, the word nonobtuse can be replaced by acute, with the modification that then each vertex projects into the relative interior of its opposite facet. In Figure 2 we illustrate this for $n=2$ and $n=3$.
\begin{center}
\includegraphics[height=5cm]{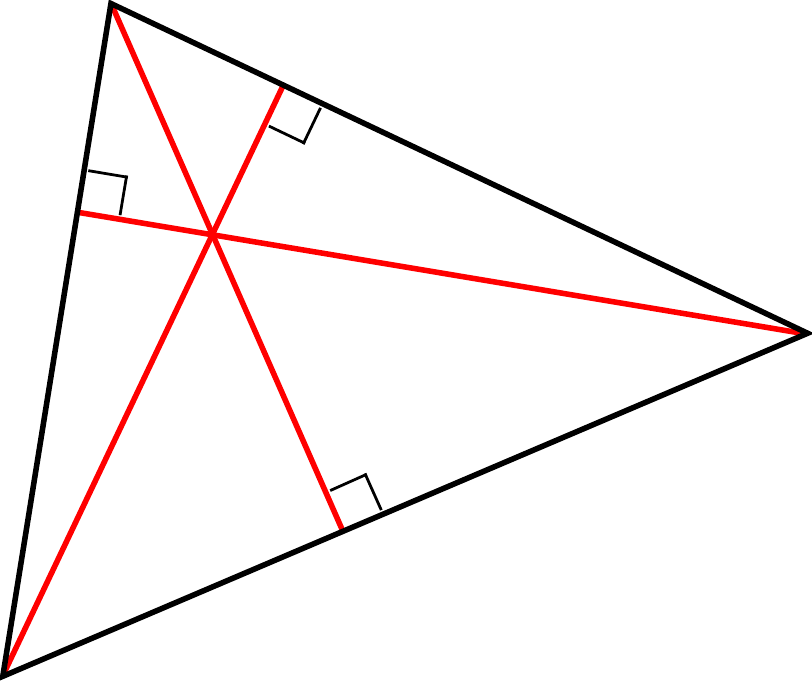}\hspace*{4mm}\includegraphics[height=5cm]{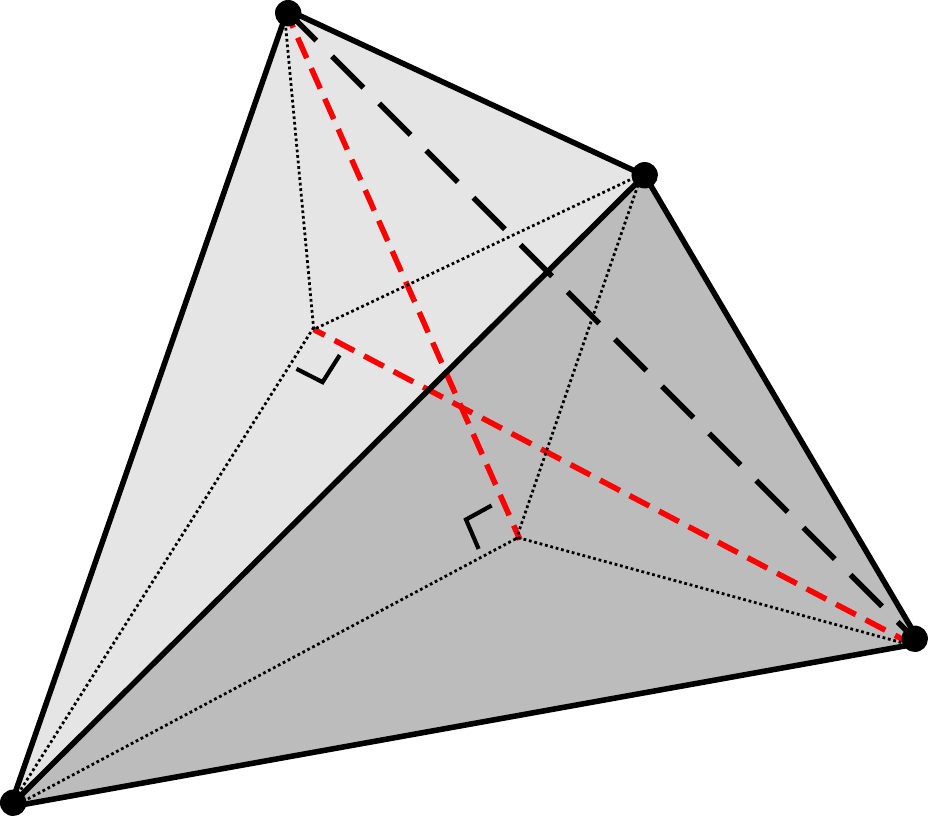}\\[2mm]
{\bf Figure 2.} Each vertex projects into the opposite facet; see Lemma \ref{pro1}.
\end{center}
Since a simplex has in total $\half n(n+1)$ dihedral angles of which at least $n$ are acute, $\half n(n-1)$ angles remain that are possibly not acute. They may, in fact, all be right. This happens if and only if $n$ of the edges of $S$ are mutualy orthogonal. Note that such an $S$ is nonobtuse.

\begin{Def}\label{def1}{\rm An $n$-simplex $S$ is called orthogonal if it has $n$ mutually orthogonal edges.}
\end{Def}
Orthogonal simplices have several special properties. Recall that a spanning tree $T$ of a graph $G$ is an acyclic connected subgraph containing all vertices of $G$. A vertex at which exactly one edge of $T$ arrives is called a leaf of the tree.

\begin{Le}\label{lem1}{\sl The orthogonal edges of an orthogonal simplex $S$ form a spanning tree $T$ for $S$.}
\end{Le}
{\bf Proof. } The vertices and edges of an orthogonal simplex $S$ form the complete graph $K_{n+1}$. The subgraph $G$ consisting of the $n$ orthogonal edges and their endpoints cannot contain a cycle due to their mutual orthogonality, and $G$ must be connected because $S$ has only $n+1$ vertices. Hence, $G$ is a spanning tree for $K_{n+1}$. \hfill $\Box$

\begin{Def}{\rm For brevity, we will call the spanning tree of orthogonal edges of an orthogonal simplex $S$ its orthogonal tree.}
\end{Def}

\subsection{Classification of orthogonal binary simplices by their orthogonal trees}\label{sect2.2}
If $S\in\BB^n$ has an exterior facet $F$, by which we mean that $F\subset \partial I^n$, the vertex opposite $F$ can be any of the $2^{n-1}$ vertices in the opposite facet of $I^n$. However, if $S$ is nonobtuse, only $n$ candidates remain. This, like more results to come, is intuitively clear, but we will provide proofs nonetheless.

\begin{Le}\label{lem3}{\sl Let $S\in\BB^n$ be nonobtuse. Then a facet $F$ of $S$ is exterior if and only if the vertex of $S$ opposite $F$ is connected to $F$ via an edge of $I^n$ orthogonal to $F$.}
\end{Le}
{\bf Proof.} Suppose that $S$ has an exterior facet $F$, thus $S$ has $n$ vertices in a facet $C$ of $I^n$. Then the remaining vertex $p$ of $S$ must lie in the facet $\ol{C}$ of $I^n$ parallel to $C$, or $S$ would be degenerate. Since each vertex $\overline{c}$ of $\ol{C}$ is connected by an edge $e$ of $I^n$ to a vertex $c$ of $C$ and $e$ is orthogonal to both $C$ and $\ol{C}$, the second item in Lemma \ref{pro1} can only hold if $p$ is connected to a vertex of $F$ by such an edge of $I^n$. Conversely, suppose that an edge $e$ of $I^n$ is an edge of $S$ and that $S$ has a facet $F$ orthogonal to $e$. Then $F$ must lie entirely in one of the two cube facets orthogonal to $e$, and is thus exterior. \hfill $\Box$\\[2mm]
The statement of this lemma is illustrated in Figure 3 below. In order for the triangle to be an exterior facet of a nonobtuse binary tetrahedron, its fourth vertex must be one of the three indicated vertices in the cube facet parallel to it.
\begin{center}
\includegraphics[width=5cm]{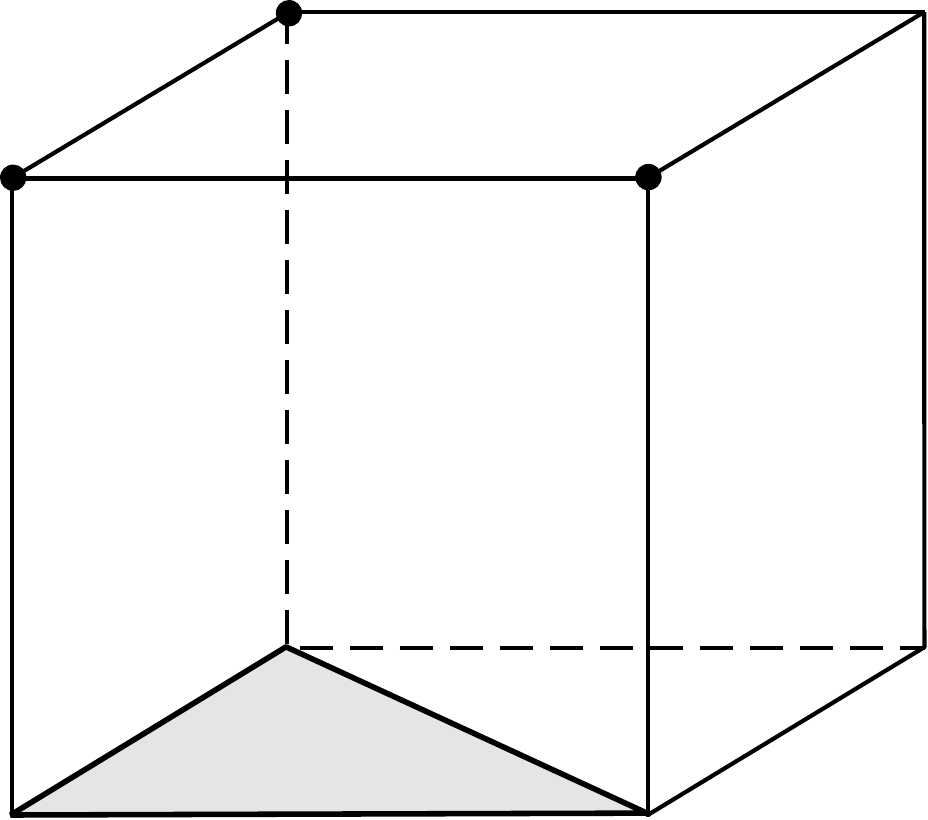}\\[2mm]
{\bf Figure 3.} Options for the remaining vertex of a nonobtuse $S\in\BB^n$ with exterior facet.
\end{center}

\begin{Le}\label{lem2}{\sl Let $S\in\BB^n$ have orthogonal tree $T$. Then each edge of $T$ that arrives at a leaf of $T$ is an edge of $I^n$.}
\end{Le}
{\bf Proof. } Let $v$ be a leaf of $T$. Without loss of generality, assume that the edge $e$ from $v$ arrives at the origin $0\in\RR^n$. Then the facet $F$ of $S$ opposite $v$ is orthogonal to $e$ and $0\in F$. As a result, each vertex of $F$ is a binary vector $q$ satisfying $\langle q,e\rangle = 0$ and must  therefore have entries equal to zero at least at those positions, where $e$ has entries equal to one. The subset of all cube vertices $p$ with two or more fixed entries equal to zero is, however, contained in an $(n-2)$-dimensional facet of $I^n$, whereas $F$ is $(n-1)$-dimensional. We conclude that $e$ cannot have two or more entries equal to one, and since $e\not=0$, it is an edge of $I^n$.\hfill $\Box$

\begin{Th}\label{co2}{\sl Let $S\in\BB^n$ be orthogonal. Then,\\[2mm]
$\bullet$ the orthogonal tree $T$ of $S$ consists of edges of $I^n$ only;\\[2mm]
$\bullet$ $S$ has $t$ exterior facets $F_1,\dots,F_t$, where $2\leq t\leq n$ is the number of leaves of $T$.\\[2mm]
Moreover, for each $j\in\{1,\dots,t\}$,\\[2mm]
$\bullet$ $F_j$ is an orthogonal $(n-1)$-simplex;\\[2mm]
$\bullet$ the orthogonal tree of $F_j$ is the subtree of $T$ induced by the vertices of $F_j$.}
\end{Th}
{\bf Proof.} By Lemma \ref{lem1} we know that $S$ has an orthogonal tree $T$. Lemma \ref{lem2} shows that each leaf $v$ of $T$ is connected to its opposite facet $F$ by an edge of $I^n$, and this edge is orthogonal to $F$. By Lemma \ref{lem3}, $F$ is exterior. This proves that $S$ has at least $t$ exterior facets. Because $F$ contains  $n-1$ of the orthogonal edges of $T$, it is itself an orthogonal $(n-1)$-simplex. Since it is contained in a facet of $I^n$, an induction argument gives that not only the edges that end in a leaf, but even all edges of $T$ are edges of $I^n$. Conversely, Lemma \ref{lem3} also shows that if $S$ has an exterior facet $F$, the vertex $p$ opposite $F$ is orthogonally connected to $F$ by an edge of the cube, and is therefore a leaf of $T$. This shows that $S$ has no more than $t$ exterior facets. The inequality $2\leq t\leq n$ expresses that each tree with $n+1$ vertices has at least $2$ and at most $n$ leaves. \hfill $\Box$\\[2mm]
Theorem \ref{co2} classifies orthogonal binary simplices according to the structure of their orthogonal tree. In Figure 4, this classification is depicted for $n\leq 6$. An arrow labeled by a number indicated how many exterior facets of a given type a certain orthogonal binary simplex has.\\[2mm]
\begin{center}
\includegraphics[width=15cm]{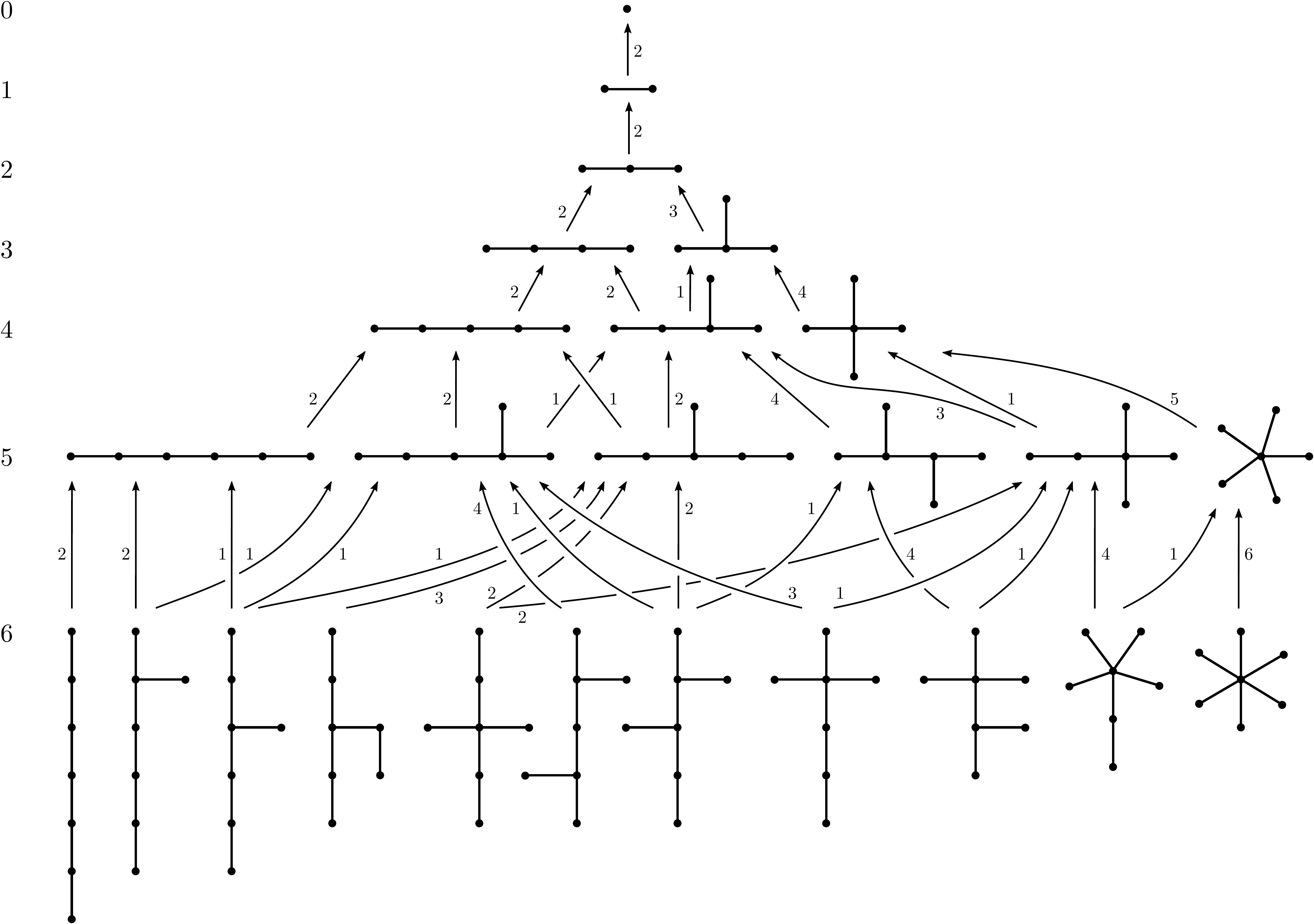}\\[2mm]
{\bf Figure 4.} Classification of orthogonal binary simplices according to their orthogonal trees.
\end{center}
For instance, for $n=4$ we see three types of orthogonal simplices. The leftmost has a path of four mutually orthogonal edges, whereas the rightmost has all its orthogonal edges meet at the same vertex. The one in the middle has two exterior facets of dimension three that have a path of three mutually orthogonal edges (see the third picture in Figure 1), and one exterior facet that has all three orthogonal edges meet at one vertex (see the rightmost picture in Figure 1). It is exactly this orthogonal simplex that will play a crucial role in what is to follow. Therefore, in the next section we pay attention to some special orthogonal simplices.

\subsection{Path simplices, cube corners, and their fake counterparts} \label{sect2.3}
A binary orthogonal simplex whose orthogonal tree is a path is called a path simplex. It has precisely two exterior facets, both of which are path $(n-1)$-simplices. A binary orthogonal simplex with $n$ exterior facets is called a cube corner. All its exterior facets are $(n-1)$-cube corners. See Figure 1 for a $3$-cube corner and a path-tetrahedron.

\begin{Def}{\rm If an orthogonal binary simplex $S$ is not a path simplex, but has an exterior facet that is a path simplex of dimension $n-1$, we call it a fake path simplex. Similarly, if $S$ is not a cube corner, but has an exterior facet that is a cube corner of dimension $n-1$, it is called a fake cube corner.}
\end{Def}

\begin{rem}{\rm Observe that by definition, only cube corners and fake cube corners have exterior facets that are $(n-1)$-cube corners, and only path simplices and fake path simplices have exterior facets that are $(n-1)$-path simplices. This will be used in Section 4.}
\end{rem}
The following two lemmas can be well understood by studying Figure 4.

\begin{Le}\label{lem8}{\sl Let $n\geq 4$. Each fake cube corner $S\in\BB^n$ has one exterior $(n-1)$-cube corner facet and $n-2$ exterior fake $(n-1)$-cube corner facets}.
\end{Le}
{\bf Proof. } Each fake cube corner $S\in\BB^n$ is the convex hull of its exterior $(n-1)$-cube corner facet $F$ and a vertex $p$ connected by an edge of $I^n$ to one of the $n-1$ leaves of the orthogonal tree of $F$. The exterior facet opposite any of the other leaves is a fake $(n-1)$-cube corner. \hfill $\Box$

\begin{Le}\label{lem9}{\sl Let $n\geq 4$. Each fake path simplex $S\in\BB^n$ has three exterior facets, of which one or two are not a path $(n-1)$-simplex}.
\end{Le}
{\bf Proof. } The orthogonal tree $T$ of $S$ contains a path between $n$ vertices $v_1,\dots,v_n$. Because $S$ itself is not a path simplex, there exists a $j\in\{2,\dots,n-1\}$ such that $v_j$ has degree $3$. If $j\geq 3$ then the exterior facet opposite $v_1$ is not a path $(n-1)$-simplex, and if $j\leq n-2$ then the exterior facet opposite $v_n$ is not a path $(n-1)$-simplex. The facet opposite the leaf connected to $v_j$ is always a path $(n-1)$-simplex.\hfill $\Box$\\[2mm]
In dimension two, as visible in Figure 4, cube corners and path simplices coincide, and thus a $3$-cube corner is a fake path simplex, and a path $3$-simplex is a fake cube corner. This peculiarity is not true anymore for dimensions $n>3$. In a sense it is the cause for the bifurcation into two distinct families of nonobtuse binary triangulations of $I^n$ for $n\geq 3$: even though each facet of $I^3$ can only be triangulated in one way, there are two distinct nonobtuse triangulations of $I^3$.

\begin{rem}\label{rem10}{\rm The middle orthogonal tree in Figure 4 for $n=4$ shows that fake cube corners and fake path simplices coincide in that dimension. We will call a simplex of this type a snake simplex. It plays an essential role in the proof of Theorem \ref{th3}.}
\end{rem}

\section{The triangulation $\TT_\NN^n$}\label{sect3}
A convex polytope $P$ with vertex $p$ can be triangulated by first triangulating all facets of $P$ that do not have $p$ as a vertex using $(n-1)$-simplices, and then taking the convex hulls of each of these simplices with $p$. This process is known \cite{Bohm} as coning off the polytope $P$ towards the vertex $p$. Here we will apply this process to $I^n$ with the cube corner opposite $p$ removed.\\[2mm]
{\bf Notation. } Let $e^n\in I^n$ be the all-ones vector. For each $j\in\{1,\dots,n\}$ write $C_j$ for the facet of $I^n$ that does not contain $e^n$ and is orthogonal to the canonical basis vector $e^n_j$ of $\RR^n$. In other words, $C_j$ is the convex hull of the binary vectors whose $j$-th coordinate equals zero. Write $K^n$ for the cube corner of $I^n$ having facets in $C_1,\dots,C_n$ and $A^n$ for the convex hull of the interior facet of $K^n$ and $e^n$. This nonobtuse simplex $A^n$ is called the antipodal of $K^n$ \cite{CoBr}.

\begin{Th}\label{th1}{\sl There exists a nonobtuse binary triangulation $\TT_\NN^n$ of $I^n$ consisting of $\NN(n)$ simplices, where
\be \NN(n) = n\NN(n-1)- n+2 \hdrie\hdrie \mbox{\rm with } \hdrie \NN(1)=1. \ee
Moreover, $K^n\in\TT_\NN^n$, and any simplex from $\TT_\NN^n$ other than $K^n$ has $e^n$ as a vertex.}
\end{Th}
{\bf Proof. } Because $I^2\setminus K^2$ is a nonobtuse triangle having $e^2$ as vertex, the statement holds for $n=2$. Assume now as inductive hypothesis that $I^{n-1}\setminus K^{n-1}$ can be triangulated into $k=\NN(n-1)-1$ nonobtuse binary simplices that all have $e^{n-1}$ as a vertex. We will show that $I^n\setminus K^n$ can be triangulated into $nk+1$ nonobtuse binary simplices that all have $e^n$ as a vertex. For this, observe that for each $j\in\{1,\dots,n\}$, the facet of $K^n$ in $C_j$ is a copy of $K^{n-1}$. Thus, among the $n+1$ facets of $I^n\setminus K^n$ that do not have $e^n$ as a vertex, there are $n$ copies of $I^{n-1}\setminus K^{n-1}$. Using the inductive hypothesis, these facets can be triangulated with $k$ nonobtuse simplices that have $e^n-e^n_j$ as vertex. Coning off $I^n\setminus K^n$ to $e^n$ results in $A^n$ and $nk$ other nonobtuse simplices. Hence, $I^n$ can be triangulated with $nk+2$ nonobtuse binary simplices, of which all but $K^n$ have $e^n$ as vertex. Since $k = \NN(n-1)-1$, we find that
\[ \NN(n)= nk+2 = n(\NN(n-1)-1)+2=  n\NN(n-1)- n+2.\]
This proves the statement. \hfill $\Box$\\[2mm]
The above proof is illustrated in Figure 5. The case $n=3$ is constructed by triangulating the three facets of $I^3$ containing the origin by copies of the triangulation for $n=2$, removing the cube corner, and coning off the resulting shape to the all-ones vertex $e^3$.
\begin{center}
\includegraphics[height=5cm]{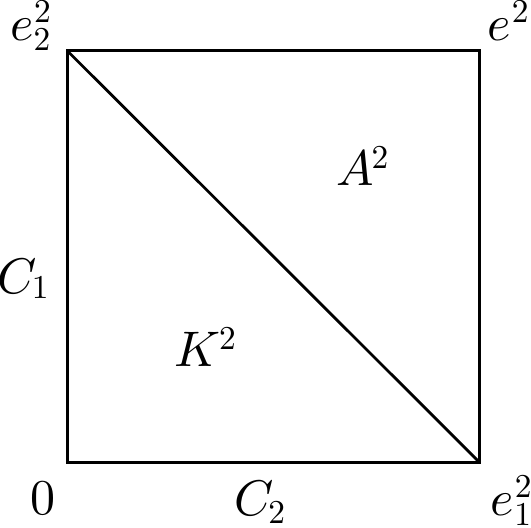}\hspace*{4mm}\includegraphics[height=5.5cm]{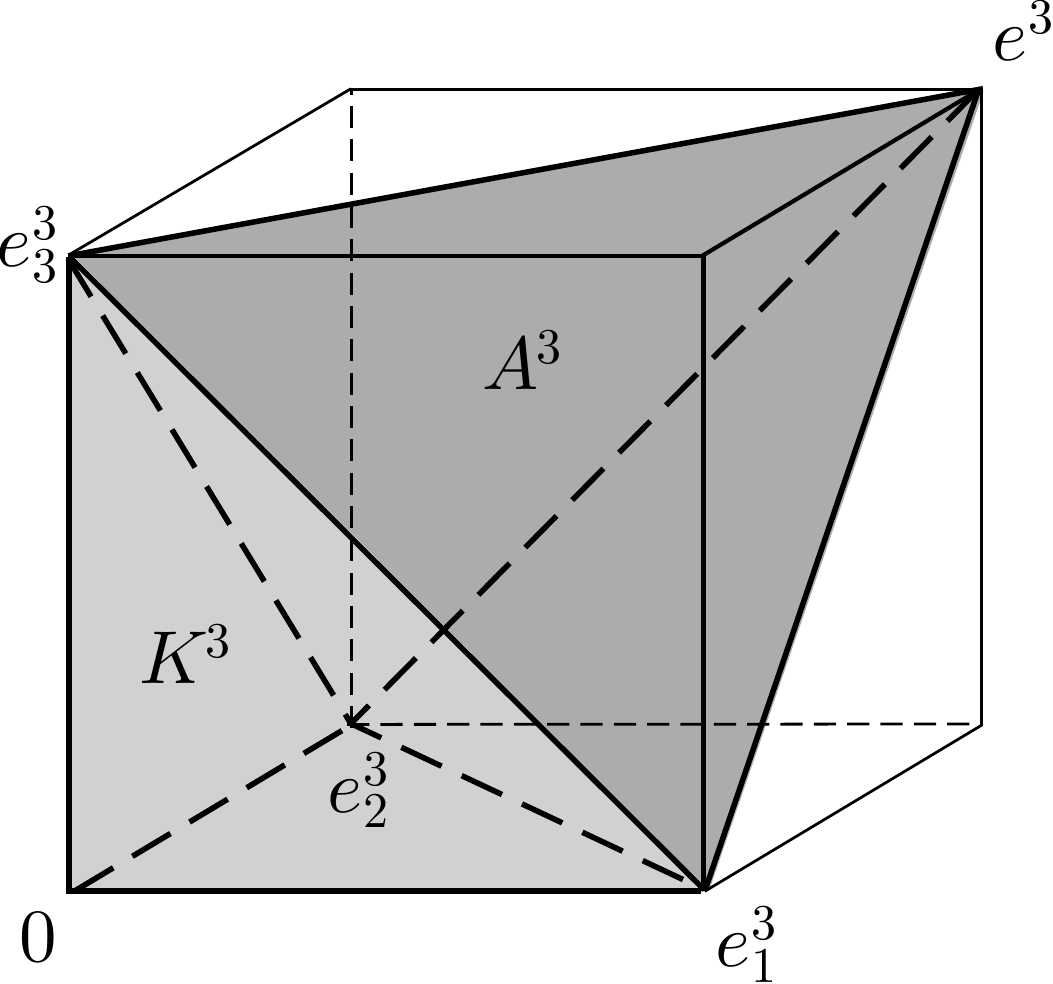}\\
{\bf Figure 5.} Illustration of the induction step in the proof of Theorem \ref{th1}.
\end{center}
We will now take a closer look at this new family of triangulations $\TT_\NN^n$. First of all, note that $\TT_\NN^n$ consists of the cube corner $K^n$; its antipodal $A^n$; $n$ simplices that are the convex hull of $e_n$ and an exterior facet of type $A^{n-1}$ which we will call of type $A^n_1$; $n(n-1)$ simplices that are the convex hull of $e^n$, one of the vectors $e^n-e^n_j$ and an $(n-2)$-dimensional facet of type $A^{n-2}$, which we shall call of type $A^n_2$, etc.\\[2mm]
For example, the triangulation $\TT_\NN^3$ consists of the cube corner $K^3$, its antipodal $A^3$, and three additional tetrahedra formed in the process of coning off towards $e^3$. These three tetrahedra are of type $A^3_1$, each being the convex hull of an exterior facet of type $A^2$ and an additional edge of length one added to its top. Similarly, $\TT_\NN^4$ consists of $K^4,A^4$, four simplices of the form $A^4_1$, each being the convex hull of an exterior facet of the form $A^3$ with a cube edge added to its top, and $4\times 3$ simplices of the form $A^4_2$, each being the convex hull of an exterior facet of the form $A^3_1$ with an additional cube edge.\\[2mm]
If finally, we write $A^n_0$ for the antipodal $A^n$, the following result is not difficult to understand.

\begin{Th}{\sl The triangulation $\TT_\NN^n$ consists of a cube corner $K^n$ together with
\be \frac{n!}{(n-k)!} \hdrie\mbox{\sl simplices of type $A^n_{k}$} \ee
for each $k\in\{0,\dots,n-2\}$. Consequently,
\be \NN(n) = 1+ n!\sum_{k=2}^{n} \frac{1}{k!} \und  \lim_{n\rightarrow\infty} \frac{\NN(n)}{n!} = {\rm e}-2. \ee
In fact, $\NN(n)$ equals the smallest integer larger than $n!({\rm e}-2)$.}
\end{Th}
{\bf Proof. } It can be verified that the first statement is a consequence of the process described in the proof of Theorem \ref{th1}. The expression for $\NN(n)$ follows from summation of over all types of simplices and their amounts, and the limit is derived from the standard exponential power series. Finally,
\[ \NN(n)-n!({\rm e}-2) = 1-n! \sum_{k=n+1}^{\infty} \frac{1}{k!} \]
and
\[ n! \sum_{k=n+1}^{\infty} \frac{1}{k!} = \frac{n!}{(n+1)!} + \frac{n!}{(n+2)!} + \dots = \frac{1}{n+1}+\frac{1}{(n+1)(n+2)}+\cdots < \sum_{j=1}^{\infty} \frac{1}{(n+1)^j} \leq 1 \]
which also proves the last statement. \hfill $\Box$

\begin{rem}{\rm Note that by construction, $\TT_\NN^n$ induces the standard triangulation on each of the $n$ facets of $I^n$ having $e^n$ as a vertex. This can be proved inductively, but also in an alternative way as will be done in Section 4.}
\end{rem}

\section{Uniqueness}\label{sect4}
In this section we will show that modulo the $n!2^n$ symmetries of the $n$-cube, the two families $\TT_\SS^n$ and $\TT_\NN^n$ of nonobtuse triangulations of $I^n$ are the only nonobtuse triangulations of $I^n$. For a given nonobtuse triangulation $\TT$ we will first show that:\\[2mm]
$\bullet$ if $\TT$ contains a path simplex, then $\TT=\TT_\SS^n$; this will be done in Section \ref{sect4.1},\\[2mm]
$\bullet$ if $\TT$ contains a cube corner, then $\TT = \TT_\NN^n$; this will be done in Section \ref{sect4.3}.\\[2mm]
In Section \ref{sect4.2}, a result is proved concerning feasible nonobtuse neighbors of a given nonobtuse binary simplex. This result will be used in Section \ref{sect4.3}, but also in Section \ref{sect4.4}, where we will show that $\TT$ contains either a path simplex or a cube corner, finishing the uniqueness proof.

\subsection{Nonobtuse triangulations $\TT$ containing a path simplex}\label{sect4.1}
By a long diagonal in a triangulation of $I^n$ we mean an edge between two antipodal vertices of $I^n$. A binary path simplex has a long diagonal of $I^n$ as an edge. In fact, also the converse is true. If a nonobtuse binary simplex contains a long diagonal, it is a path simplex.

\begin{Le}\label{lem4}{\sl If a nonobtuse $S\in\BB^n$ has a long diagonal of $I^n$ as an edge, $S$ is a path simplex.}
\end{Le}
{\bf Proof. } Without loss of generality, assume that $S\in\BB^n$ is nonobtuse with both the origin and $e^n$ as vertices. The facet $F$ of $S$ opposite $e^n$ is contained in a hyperplane $H$ through the origin. Let $q$ with $\|q\|=1$ be the vector normal to $H$ with orientation determined by $q^\top e^n>0$. Note that $q^\top e^n\not=0$ because $e^n\not\in H$. Now, consider the orthogonal projection $p=e^n - q(q^\top e^n)$ of $e^n$ on $H$.
\begin{center}
\includegraphics[height=6cm]{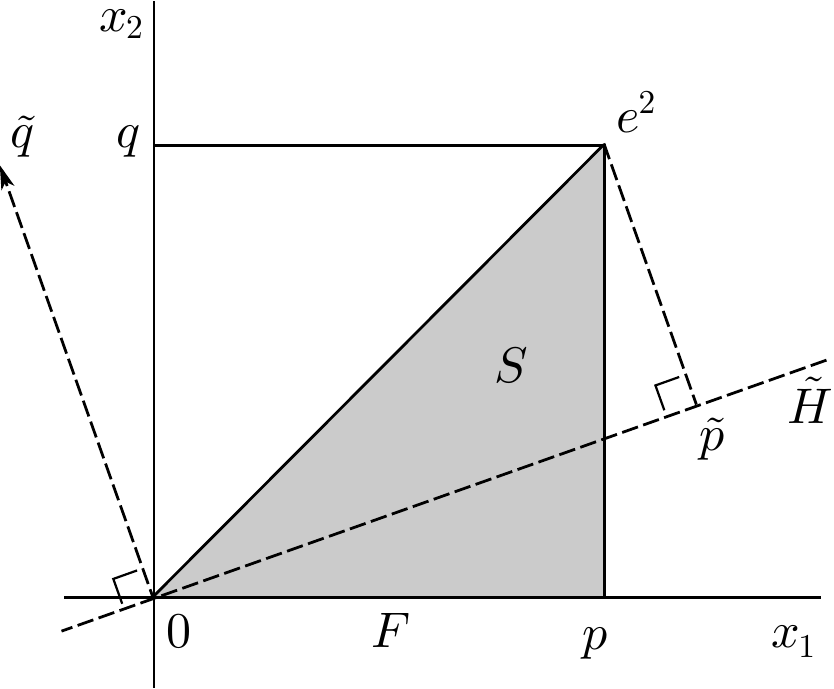}\\
{\bf Figure 6.} The projection $\tilde{p}$ of $e$ on $H$ does not lie in $I^n$.
\end{center}
We will first show that if $q$ is not a standard unit basis vector, then $p\not\in F$ and $S$ is obtuse by Lemma \ref{pro1}. For this, write $q_j=(e_j^n)^\top q$ and $p_j=(e_j^n)^\top p$. Suppose that $j$ is such that $q_j<0$. Then $p_j = 1-q_j(q^\top e_n) > 1$ due to $q^\top e^n>0$,  and thus $p$ is not in $I^n$, and also not in $F$, which contradicts Lemma \ref{pro1}. This situation is indicated in Figure 6 by $\tilde{q},\tilde{p}$ and $\tilde{H}$. We conclude that a $j$ with $q_j<0$ does not exist. Hence, $q$ is nonnegative. Then we can use an argument from the proof of Lemma \ref{lem2}. Indeed, each vertex of $F$, being a binary vector orthogonal to $q$, must have a zero entry at each position where $q$ has a positive entry. Since $F$ has dimension $n-1$, this is only possible if $q$ has exactly one nonzero entry. We conclude that $q=e_j^n$ for some $j\in\{1,\dots,n\}$. Lemma \ref{lem3} then shows that $F$ is an exterior facet. Since $p\in F$, we have that $F$ contains the long diagonal of the cube facet that contains $F$. Since $F\in\BB^{n-1}$, the proof can now easily be finished by induction.\hfill $\Box$

\begin{Th}\label{co3}{\sl $\TT_\SS^n$ is the only nonobtuse binary triangulation of $I^n$ containing a path simplex}.
\end{Th}
{\bf Proof. } Let $\TT$ be a nonobtuse binary triangulation of $I^n$ containing a path simplex $S\in\BB^n$. Each interior facet of $S$ contains the same long diagonal. According to Lemma \ref{lem4}, these interior facets are met by other path simplices, leading to the standard triangulation into $n!$ path simplices. Hence, $\TT=\TT_\SS^n$. \hfill $\Box$

\subsection{Nonobtuse neighbors of interior facets}\label{sect4.2}
In this section we will show that an interior simplicial facet $F$ whose normal has no zero entries can have at most two nonobtuse neighbors, one on each side of $F$. To do so, define for given $w\in\RR^n$ the translation
\be \tau_w: \RR^n\rightarrow \RR^n: \hdrie x \mapsto x-w,\ee
that maps $w$ onto the origin. For each orthant $E$ of $\RR^n$ there exists exactly one vertex $v$ of $I^n$ such that $\tau_v(I^n)\subset E$, and if $v_1$ and $v_2$ are opposite vertices then $\tau_{v_1}(I^n)$ and $\tau_{v_2}(I^n)$ lie in opposite orthants. Explicitly, it can be easily verified that if $v$ is a vertex of $I^n$, then $\tau_v$ maps $I^n$ from the nonnegative orthant into the orthant defined by the relations
\be\label{new1} x_j\geq 0 \Leftrightarrow v^\top e_j^n = 0 \und x_j\leq 0 \Leftrightarrow v^\top e_j^n=1, \hdrie j\in\{0,\dots,n\}.\ee
For instance, if $w=(1,1,0)$ then $\tau_w$ maps $I^3$ into the octant $x_1\leq 0, x_2 \leq 0, x_3 \geq 0$.\\[2mm]
From this observation we derive a useful lemma, in which the following notation is used. For $u,v\in\RR^n$, we write 
\be \ell_u(v)=\{v+\lambda u\sth \lambda\in\RR\} \ee
for the line through $v$ parallel to the vector $u$.

\begin{Le}\label{newlem1}{\sl If $u$ has no zero entries then $\ell_u(v)\cap I^n\not=\{v\}$  for exactly two antipodal vertices $v=v_1$ and $v=v_2$ of $I^n$. Moreover, the entries of $v_1$ are determined by}
\be\label{new2}  v_1^\top e_j^n=0 \Leftrightarrow u^\top e_j^n>0 \und v_1^\top e_j^n = 1\Leftrightarrow u^\top e_j^n<0, \hdrie j\in\{1,\dots,n\}, \ee
and $v_2=e^n-v_1$ is the antipodal of $v_1$.
\end{Le}
{\bf Proof. } If $u\in\RR^n$ has no zero entries, it is an element in the interior of an orthant $E_u$ of $\RR^n$, and $\ell_u(0)$ is a line through the origin. This line only intersects the interiors of $E_u$ and its opposite orthant $E_{-u}$. Now, for a given vertex $v$ of $I^n$ we have that
\be \ell_u(0) = \tau_v(\ell_u(v)).\ee
Therefore, $\ell_u(v)$ intersects $I^n$ in points other than $v$ if and only if $\tau_v$ maps $I^n$ into either $E_u$ of $E_{-u}$. The relations (\ref{new2}) follow immediately from (\ref{new1}).  \hfill $\Box$

\begin{Co}\label{newcor1}{\sl Let $F$ be an interior facet of a simplex $S\in\BB^n$ having a normal vector $\nu$ without zero entries. Then there exist at most two nonobtuse simplices $S_1$ and $S_2$ that have $F$ as common facet, and the vertex of $S_1$ that does not belong to $F$ is antipodal to the vertex of $S_2$ that does not belong to $F$.}
\end{Co}
{\bf Proof. } Let $v\not\in F$ be a vertex of $I^n$. A necessary condition for $v$ to project orthogonally onto $F$ is that $\ell_\nu(v)\cap I^n$ consists of more than only the vertex $v$ itself. Lemma \ref{newlem1} shows that this condition is satisfied by at most two antipodal vertices of $I^n$.\hfill $\Box$\\[2mm]
Thus, in words, if the normal to an interior facet $F$ has no zero entries, $F$ has at most two nonobtuse neighbours, one on each side of $F$.
 
\begin{Le}\label{lem5}{\sl For any nonobtuse binary triangulation $\TT$ of $I^n$ we have}
\[ K^n\in\TT \Leftrightarrow A^n\in\TT. \]
\end{Le}
{\bf Proof. } This follows immediately from Corollary \ref{newcor1}, because the normal to the interior facet of $K^n$ is a multiple of $e^n$. \hfill $\Box$

\subsection{The cube corner triangulation}\label{sect4.3}
In Section \ref{sect4.1} we saw that the presence of one path simplex in a nonobtuse triangulation $\TT$ of $I^n$ implies that $\TT=\TT_\SS^n$. Here we prove that if $\TT$ contains a cube corner, then $\TT=\TT_\NN^n$.\\[2mm]
Roughly speaking, the proof shows that if one cube facet $C_n$ is triangulated with $\TT_\NN^{n-1}$, then any additional edge sprouting from the vertex of $C_n$ furthest away from the origin, must go vertically up in the new dimension. This, however, is precisely what happens in the process of coning off to $e_n^n$, and thus leads to $\TT_\NN^{n}$.

\begin{Th}\label{th2}{\sl $\TT_\NN^n$ is the only nonobtuse binary triangulation of $I^n$ containing a cube corner}.
\end{Th}
{\bf Proof. } The statement holds trivially for $n=2$. As inductive hypothesis, assume that the only nonobtuse triangulation of $I^{n-1}\setminus K^{n-1}$ is $\TT_\NN^{n-1}\setminus \{K^{n-1}\}$. Let $\TT$ be a nonobtuse triangulation of $I^n\setminus K^n$. Consider, without loss of generality, the bottom facet $C_n$ of $I_n$. By the induction hypothesis, its triangulation is a copy of $\TT_\NN^{n-1}$. Let $S^{n-1}$ be any simplex in $C_n$ from $\TT_\NN^{n-1}\setminus \{K^{n-1}\}$. Then $e^n-e^n_n$ is a vertex of $S^{n-1}$, as is shown in Theorem \ref{th1}. Naturally, $S^{n-1}$ is the exterior facet of a simplex $S^n\in\TT$. Write $v$ for the vertex of $S^n$ opposite $S^{n-1}$. Below we will demonstrate that $v=e^n$. This will finish the proof, because then $S^n$ is obtained by coning off $S^{n-1}$ to $e^n$ and thus a simplex from $\TT_\NN^n$, and $S^{n-1}$ was arbitrary. To start, first observe that due to Lemma \ref{pro1}, $v\not=e_n^n$ because the projection of $e_n^n$ on $C_n$ is the origin, which is not a vertex of $S^{n-1}$. Thus, $v$ is a vertex in a facet $C$ of $I^n$ that contains $e^n-e_n^n$. In particular, the edge $\ell$ of $S^n$ between $v$ and $e^n-e^n_n$ is also contained in $C$. Since obviously, $C\not=C_n$, we assume without loss of generality that $C$ is the facet of $I^n$ contained in the plane $x_1=1$. Since $K^n\in\TT$, by Lemma \ref{lem5} also $A^n\in\TT$, and thus the long diagonal of $C$ between $e^n_1$ and $e_n^n$ is an edge in $\TT$. Lemma \ref{lem4} and Theorem \ref{co3} yield that $C$ is triangulated by path simplices. But in the triangulation by path simplices, the only edge from $e^n-e_n^n$ to a vertex in the cube facet parallel to $C_n$ is the cube edge that connects them. Therefore, we conclude that $v=e^n$.\hfill $\Box$

\subsection{Nonobtuse triangulations without path simplices and cube corners}\label{sect4.4}
We will prove that the triangulations in the title of this section do not exist. This will be done separately for $n=4$ and $n=5$ and then for $n\geq 6$ by induction. For $n=4$, the proof is relatively simple and based on a counting argument of the visible exterior facets. For $n=5$, we need a rather involved argument based on a double application of Corollary \ref{newcor1} in order to show that an alternative triangulation fails due to lack of nonobtuse neighbours. Fortunately, the induction step is then relatively straightforward again.

\begin{Th}\label{th3}{\sl The only nonobtuse binary triangulations of $I^4$ are $\TT_\SS^4$ and $\TT_\NN^4$}.
\end{Th}
{\bf Proof. } Assume without loss of generality, based on Theorems \ref{co3} and \ref{th2}, that $\TT$ is a nonobtuse triangulation of $I^4$ without cube corners and path simplices. Lemma \ref{pro1} shows that $\TT$ induces a nonobtuse triangulation of each of the eight facets of $I^4$. There are only two possibilities for the triangulation of each facet: either in six path tetrahedra, or in four cube corners and a regular simplex. Thus, for some $p\in\{0,\dots,8\}$, we observe in total $6p$ path tetrahedra and $4(8-p)$ cube corners at the boundary of $I^4$. By Remark \ref{rem10}, these can only be the exterior facets of snake simplices. Each snake simplex has two exterior path tetrahedra and one exterior cube corner. However, the numbers $6p$ and $4(8-p)$ are not in the ratio $2:1$ for any $p\in\{0,\dots,8\}$. This shows that a triangulation $\TT$ without path simplices and without cube corners does not exist. \hfill $\Box$\\[2mm]
Unfortunately, setting up a similar counting argument for the case $n=5$ does not (seem to) lead to an immediate contradiction. We did however succeed in a number of different ways to prove Theorem \ref{th5} below. Its given proof does not rely on computer aid, but is rather lengthy. After the proof, we comment on alternatives that use the computer.

\begin{Th}\label{th5}{\sl The only nonobtuse binary triangulations of $I^5$ are $\TT_\SS^5$ and $\TT_\NN^5$}.
\end{Th}
{\bf Proof. } If a nonobtuse triangulation $\TT$ of $I^5$ induces a triangulation of each facet of $I^5$ consisting of path $4$-simplices, then Lemma \ref{lem9} shows that this is due to the presence of path simplices in $\TT$, and not due to the presence of fake path simplices. Thus, in combination with Theorem \ref{co3} this yields that $\TT=\TT_\SS^5$. Assume therefore that at least one facet of $I^5$ is triangulated differently.  According to Theorem \ref{th3}, these facets must be triangulated with copies of $\TT_\NN^4$. As a result $\TT$ must contain a nonobtuse simplex $S$ with an antipodal $4$-simplex $A^4$ as exterior facet $F_e$. According to Lemma \ref{lem3}, the vertex $p$ of $S$ opposite $A^4$ can be at no more than five different locations in the facet of $I^5$ parallel to the cube facet that contains $F_e$. First, suppose that $p$ orthogonally projects onto the top of the antipodal exterior facet. Then $S$ has a facet $F$ consisting of the convex hull of $p$ and the four vertices of $A^4$ onto which $p$ does not orthogonally project. But this is a facet of the antipodal simplex $A^5$, and a simple computation shows that its normal does not contain zero entries. By Corollary \ref{newcor1}, $A^5$ must be an element of $\TT$. By Lemma \ref{lem5} this implies the presence of a $5$-cube corner and thus by Theorem \ref{th2} we have $\TT=\TT_\NN^5$. Secondly, suppose that $p$ orthogonally projects on one of the other vertices of the exterior $4$-antipodal facet. We will show that this cannot lead to a nonobtuse triangulation $\TT$. Without loss of generality, we will consider the following explicit situation. We embed $A^4$ in the facet $x_5=0$ of $I^n$ with its top at the origin, and with $p$ above one of the vertices not equal to the origin. Then $S$ is the simplex having as vertices the origin and the $5$ columns of the matrix $P$ below, and with normal $\nu$ to its facet $F_2$ opposite the vertex in the second column of $P$, 
\[ P = \left[\begin{array}{ccccc}
0 & 0 & 1 & 1 & 1 \\ 1 & 1 & 0 & 1 & 1 \\
1 & 1 & 1 & 0 & 1\\ 1 & 1 & 1 & 1 & 0 \\ 1 & 0 & 0 & 0 & 0\end{array}\right] \hdrie\mbox{\rm with }\hdrie \nu=\left[\begin{array}{r} -2 \\ 1 \\ 1 \\ 1 \\ -3\end{array}\right] \und P^\top \nu = 3e_2.\]
Since $P^\top \nu$ has all entries zero apart from the second, $\nu$ is orthogonal to the span of all but the second column of $P$. This is a hyperplane that contains $F_2$. Thus, $\nu\perp F_2$.  Now, since $\nu$ has no zero entries, according to Corollary \ref{newcor1}, the only candidate for a nonobtuse neighbour of $S$ sharing the facet $F_2$ is the simplex $T$ whose vertices are the origin and the columns of the matrix $Q$ obtained by negating the second column of $P$,
\be \label{new3} Q = \left[\begin{array}{ccccc}
0 & 1 & 1 & 1 & 1 \\ 1 & 0 & 0 & 1 & 1 \\
1 & 0 & 1 & 0 & 1\\ 1 & 0 & 1 & 1 & 0 \\ 1 & 1 & 0 & 0 & 0\end{array}\right] \hdrie\mbox{\rm with }\hdrie \mu = \left[\begin{array}{r} 1 \\ -3 \\ 2 \\ 2 \\ -1\end{array}\right] \und Q^\top \mu = 5e_3.\ee
Recall that negating a vertex gives its antipodal. It can be verified that $T$ is nonobtuse. This is however not needed to finish the argument, because if $T$ would be obtuse, then the argument would be finished already by lack of a nonobtuse neighbour to an interior facet of $S$. The point is that the facet $G$ of $T$ opposite the vertex in the third column of $Q$ has no nonobtuse neighbour. The normal $\mu$ of this facet $G$, given explicitly in (\ref{new3}), contains no zero entries, and according to Corollary \ref{newcor1} the only candidate nonobtuse neighbour of $T$ sharing the facet $G$ is the simplex whose vertices are the origin together with the vertices in the columns of the matrix  
\[ R = \left[\begin{array}{ccccc}
0 & 1 & 0 & 1 & 1 \\ 1 & 0 & 1 & 1 & 1 \\
1 & 0 & 0 & 0 & 1\\ 1 & 0 & 0 & 1 & 0 \\ 1 & 1 & 1 & 0 & 0\end{array}\right],\]
obtained by negating the third column of $Q$. This simplex is however obtuse, because the exterior normals $\nu_3$ and $\nu_4$ between the two facets opposite its vertices in the third and fourth column make an acute angle. Indeed, 
\[ R^\top \nu_3 = -4e_3, \hdrie R^\top \nu_4=-4e_4, \hdrie \nu_3^\top \nu_4 = 1, \hdrie\mbox{\rm where} \hdrie \nu_3 = \left[\begin{array}{r} 1 \\ -3 \\ 2 \\ 2 \\ -1\end{array}\right] \und \nu_4 = \left[\begin{array}{r} -1 \\ -1 \\ 2 \\ -2 \\ 1\end{array}\right]. \]
Note that since the third entry of $R^\top \nu_3$ is negative, the vertex in the third column of $R$ is not in the same halfspace with respect to the orthogonal complement of $\nu_3$ than $\nu_3$ itself, showing that $\nu_3$ is indeed an outward pointing normal to the facet. Since dihedral angles between two facets are by definition equal to $\pi$ minus the angle between their outward normals, the simplex is obtuse. Thus, $p$ necessarily must project orthogonally onto the top of $A^4$. This proves the statement. \hfill $\Box$

\begin{rem}{\rm The proof of the above theorem can be given with the aid of the computer in various ways. For instance, we wrote a program that simply lists all nonobtuse simplices in $I^5$ and aims to triangulate $I^5$ with them. The surprising outcome (together with the corresponding results for $I^4$) is what initiated this paper. Alternatively, looking at the list of nonobtuse simplices in $I^5$, apart from the path-simplices and the antipodal simplex, there turned out to be only one other nonobtuse simplex that contained the cube's midpoint. This simplex however had interior facets without nonobtuse neighbours.}
\end{rem}  

\begin{Le}\label{lem11}{\sl Let $n\geq 5$. Then neither $\TT_\SS^n$ nor $\TT_\NN^n$ contains a fake cube corner.}
\end{Le}
{\bf Proof. } Trivially, $\TT_\SS^n$ contains no fake cube corners. Secondly, $\TT_\NN^n$ induces $n$ facets of $I^n$ to be triangulated by copies of $\TT_\SS^{n-1}$, and these do not contain $(n-1)$-cube corners. The remaining $n$ facets are each triangulated with a copy of $\TT_\NN^{n-1}$ and thus each show one exterior $(n-1)$-cube corner facet. These are however induced by the $n$-cube corner that is part of $\TT_\NN^n$. Therefore, $\TT_\NN^n$ contains no fake cube corner.\hfill $\Box$\\[2mm]
Note that $\TT_\NN^4$ does indeed contain $12$ fake cube corners, also called snakes in this paper.

\begin{Th}{\sl The only nonobtuse binary triangulations of $I^n$ for $n\geq 6$ are $\TT_\SS^n$ and $\TT_\NN^n$.}
\end{Th}
{\bf Proof. } For $n=5$ this is proved in Theorem \ref{th5}. As inductive hypothesis, assume that the statement holds for $I^n$, and let $\TT$ be a nonobtuse binary triangulation of $I^{n+1}$. Then by the hypothesis, $\TT$ only induces copies of $\TT_\SS^n$ and $\TT_\NN^n$ on the facets of $I^{n+1}$. If $\TT$ induces a copy of $\TT_\SS^n$ on each of the facets of $I^{n+1}$, then Lemma \ref{lem9} shows that this is not due to fake path simplices but due to real path simplices, and Theorem \ref{co3} yields that  $\TT=\TT_\SS^{n+1}$. On the other hand, if $\TT$ induces a copy of $\TT_\NN^n$ on at least one of the facets of $I^{n+1}$, then it contains an $(n+1)$-simplex $S$ with a cube corner facet. By definition, $S$ is either a true cube corner or a fake cube corner. If $S$ is a fake cube corner, then due to Lemma \ref{lem8}, which shows that $S$ has fake cube corner facets, $\TT$ would generate a fake cube corner in a facet of $I^{n+1}$ . However, since by the hypothesis, $\TT$ only induces copies of $\TT_\SS^n$ and $\TT_\NN^n$ on the facets of $I^{n+1}$, this contradicts Lemma \ref{lem11}. Thus, $S$ is a true cube corner, and $\TT=\TT_\NN^{n+1}$ by Theorem \ref{th2}.\hfill $\Box$\\[2mm]
Therefore, finally, we have proved the statement in the title of this paper.

\subsection*{Acknowledgments}
We would like to express our sincere gratitude for the comments and suggestions by both referees, who obviously read the paper with great care and attention for details. In particular, one of them suggested to include Lemma \ref{newlem1}, and outlined its proof. Thanks to this Lemma we were able to prove Theorem \ref{th5} without the aid of the computer, making the proof much more attractive. The referees were also very helpful in improving the presentation of the revised version of this paper. \\[2mm]
Michal K\v{r}\'{\i}\v{z}ek was supported by grant IAA 100190803 of the Grant Agency of the Academy of Sciences of the Czech Republic.
 
\subsection*{Authors' addresses}
Jan Brandts, Sander Dijkhuis, and Vincent de Haan: Korteweg-de Vries Institute for Mathematics, Faculty of Science, University of Amsterdam, Science Park 904, P.O. Box 94248, 1090 GE Amsterdam, Netherlands. E-mail: janbrandts@gmail.com, sander.dijkhuis@gmail.com, vincentdehaan@gmail.com\\[2mm]
Michal K\v{r}\'{\i}\v{z}ek, Institute of Mathematics, Academy of Sciences, \v{Z}itn\'a 25, 115 67 Prague 1, Czech Republic. E-mail: krizek@math.cas.cz\\[2mm]
Corresponding author: Jan Brandts.

\end{document}